\documentclass[12pt]{amsart}
\usepackage{amscd}

\newtheorem{thm}{Theorem}[section]

\newtheorem{lma}[thm]{Lemma}
\theoremstyle{definition}

\def \ph{\varphi}

\def \gr{graded}

\def \TV{{\mathcal T}(V)}

\def \refeq#1{equation (\ref{#1})}
\def \ra{\rightarrow}

\def \hom{\mbox{\rm Hom}}
\def \ie{\hbox{\it i.e.}}

\def \tns{\otimes}

\def \mwedge{\wedge\cdots\wedge}
\def \mplus{+\cdots+}
\def \mcom{,\cdots,}
\def \k{\mbox{$\mathfrak K$}}

\def \Z{\mbox{$\mathbb Z$}}

\def\br#1#2{[#1,#2]}

\def\zt{\mbox{$\Z_2$}}
\def\ztz{\mbox{$\Z_2\times\Z$}}

\def\sh{\operatorname{Sh}}

\def\inv{^{-1}}
\def\d{d}
\def\td{\tilde\d}
\def\id{\operatorname{Id}}
\def\im{\operatorname{Im}}
\def\A{\mbox{$\mathcal A$}}
\def\B{\mbox{$\mathcal B$}}
\def\L{L}
\def\O{\mathcal O}
\def\LA{\mbox{$\L_{\A}$}}
\def\LB{\mbox{$\L_{\B}$}}

\def\WA{\mbox{$\W_{\A}$}}

\def\m{\mbox{$\mathfrak m$}}

\def\Lm{\mbox{$\L\htns\m$}}

\def\V{V}
\def\W{W}

\def\der{\operatorname{Der}}
\def\coder{\operatorname{Coder}}
\def\ainf{\mbox{$A_\infty$}}
\def\linf{\mbox{$L_\infty$}}
\def\and{\mbox{ \rm and }}
\def\T{\mathcal T}
\def\TV{\T(V)}
\def\TW{\T(W)}

\def\s#1{(-1)^{#1}}
\def\invlim{\operatorname{\overleftarrow\lim}}
\def\tl{\tilde\lambda}
\def\HL{H(\L)}
\def\ZL{Z(\L)}
\def\BL{B(\L)}

\def\htns{\hat\tns}
\def\SW{S(W)}
\def\SWA{\SW\tns\A}
\def\htns{\hat\tns}
\def\tanA{\hbox{TA}}
\def\H{\mathcal H}
\author{Alice Fialowski}
\address{The Lor\'and E\"otv\"os University\\
Budapest, Hungary} \email{fialowsk@cs.elte.hu}
\author{Michael Penkava}
\address{University of Wisconsin\\
Eau Claire, WI 54702-4004} \email{penkavmr@uwec.edu}

\thanks{The research of the first author was partially supported by the
grants OTKA T030823, T23434, T29535 and FKFP 0170/1999, and the second
author by an NRC grant, as well as grants from the University of
Wisconsin-Eau Claire.}
\title{Deformation Theory of Infinity Algebras}
\begin{document}
\setlength{\multlinegap}{0pt}
\nocite{ps2,pen1,ls,pen2,pen3,kon,fm,mar,mar2,ksv,sta1,sta2,sta3,getz,getz2,
ge_ka1,
gers,lod,hoch,fi,fi2,ff2,ff3,pw,ger,harr,ss,bar}

\begin{abstract}
 This work explores the deformation theory of algebraic structures in a
very general setting. These structures include commutative,
associative algebras, Lie algebras, and the infinity versions of
these structures, the strongly homotopy associative and Lie
algebras. In all these cases the algebra structure is determined by an
element of a certain graded Lie algebra which plays the role of
a differential on this algebra. We work out the deformation
theory in terms of the Lie algebra of coderivations of an
appropriate coalgebra structure and construct a universal
infinitesimal deformation as well as a miniversal formal
deformation. By working at this level of generality, the main
ideas involved in deformation theory stand out more clearly.
\end{abstract}
\maketitle

\section{Introduction}
In this paper we explore the notion of deformations of algebraic
structures, in a very general setting, which is applicable to any
algebraic structure determined by a differential on a
graded Lie algebra. Examples of such structures include
associative algebras, which are determined by a differential on
the Lie algebra of coderivations of the tensor coalgebra, Lie
algebras, which are determined by a coderivation on the
coderivations of the exterior coalgebra, as well as the infinity
versions of these structures, which are determined by more
general differentials than the quadratic ones which give the
usual associative and Lie algebras.

Hochschild cohomology of associative algebras was first described
in \cite{hoch}, and used to study the infinitesimal deformations
of associative algebras in \cite{gers}, \cite{ger}. In \cite{gers} the
Gerstenhaber bracket was defined on the space of cochains of an
associative algebra, which equips the space of cochains with the
structure of a graded Lie algebra. Moreover, it was shown that
the cochain determining the associative algebra structure is a
differential on this Lie algebra, whose homology coincides with
the Hochschild cohomology of the associative algebra. 
Jim Stasheff discovered (see \cite{sta4}) that this construction could be understood more
simply by means of a natural identification of the space of
cochains of an associative algebra with the coderivations of the
tensor coalgebra on the underlying vector space. With this
identification, the Gerstenhaber bracket coincides with the
bracket of coderivations. Thus the initially mysterious existence
of a Lie algebra structure on the space of cochains of an
associative algebra was unraveled. In Stasheff's construction, the
associative algebra structure is determined by a 2-cochain, \ie,
a quadratic coderivation, which has square zero, so that it is a
codifferential on the \emph{tensor} coalgebra. Infinitesimal
deformations of the algebra are determined by quadratic cocycles
with respect to the homology determined by the codifferential, in
other words, by the second cohomology group. The connection
between deformation theory and cohomology is completely
transparent in this framework. It also turned out (see \cite{sta1})
 that replacing
the quadratic codifferential by a more general codifferential
determines an interesting algebraic structure, which is called an
\ainf\ algebra, or strongly homotopy associative algebra. \ainf\
algebras first were described in \cite{sta1}, \cite{sta2}, and have
appeared in both mathematics and physics. (See \cite{getz}, \cite{kon},
\cite{mar}, \cite{ps2}, \cite{ksv}, \cite{ge_ka1}.)

In a parallel manner, the notion of a Lie algebra can be described in
terms of the coderivations of the \emph{exterior} coalgebra of a vector
space, with the Lie algebra structure being given by a quadratic
codifferential, and the Eilenberg-Chevalley cohomology of the Lie
algebra being given by the homology determined by this codifferential.
In addition, the homology of a Lie algebra inherits the structure of a
graded Lie algebra from the bracket of coderivations. Lie algebras
generalize to \linf\ algebras (strongly homotopy Lie algebras). They
first appeared in \cite{ss}, and also have applications in both
mathematics and physics. (See \cite{ls}, \cite{sta3}, \cite{lm},
\cite{pen1}, \cite{pen2}, \cite{pen3}, \cite{bar}.)

 The main feature of
this description is that in the cases described above, the algebraic
structure is determined by an element of a certain graded Lie algebra,
which plays  the role of the differential on this algebra, and the
cohomology of the algebra is simply the homology of this differential.
Infinitesimal deformations of the algebraic structure are completely
determined by the cohomology, which has the structure of a graded Lie
algebra. There are many other algebraic structures which fit this basic
framework. For example, commutative algebras are determined by a
quadratic codifferential in the space of coderivations of the Lie
coalgebra associated to the tensor coalgebra. Similarly, deformations
of associative or Lie algebras preserving an invariant inner product,
are determined by cyclic cohomology, which is given by a differential
graded Lie algebra on a space of cyclic cochains (see \cite{lod},
\cite{pen1},\cite{mar2}).

The notion of deformations of a Lie algebra with base given by a
commutative algebra was described in \cite{fi} and used in
\cite{ff3} in order to study some examples of singular
deformations (see also \cite{fi2}). In \cite{ff2} a universal infinitesimal
deformation of a Lie algebra was constructed, as well as a
miniversal formal deformation.Our purpose in this article is to
 generalize these
two constructions to the general setting of an algebraic
structure which is determined by a differential on a graded Lie
algebra, as in the structures described above.

First, we will give a description of deformation theory in terms
of the Lie algebra of coderivations of an appropriate coalgebra
structure. In this general setting we will construct a universal
infinitesimal deformation, as well as a miniversal formal
deformation. By working at this level of generality,  the main
ideas involved in the construction stand out more clearly.

\section{Infinity algebras}
As a preliminary exercise, we first translate the notions of Lie
algebras and associative algebras into descriptions in terms of
the language of codifferentials on symmetric and tensor algebras,
which allows us to give simple definitions of their
generalizations into \linf\ and \ainf\ algebras, as well as
making it possible to describe deformation theory in a uniform
manner.

Let \k\ be a field and \V\ be a \zt-graded \k-vector space. (Some
may prefer a \Z-grading on \V.) Then by the exterior algebra
$\bigwedge V$ we mean the quotient of the (restricted) tensor
algebra $\TV=\displaystyle\bigoplus_{n=1}^\infty V^n$ by the
graded ideal generated by elements of the form $u\tns v +\s{uv}
v\tns u$ for homogeneous elements $u$, $v\in V$ (where $\s{uv}$ is
an abbreviation for $\s{|u||v|}$, with $|u|$ denoting the parity
of $u$). The exterior algebra may also be considered as the
symmetric algebra arising from a natural \ztz-grading on $\TV$,
and it is \ztz-graded commutative in this sense.

If $\sigma$ is a permutation in $\Sigma_n$, then
\begin{equation*}
v_{\sigma(1)}\mwedge v_{\sigma(n)}= \s{\sigma}\epsilon(\sigma)
v_1\mwedge v_n
\end{equation*}
where $\s{\sigma}$ is the sign of $\sigma$, and
$\epsilon(\sigma)$ is a sign depending on both $\sigma$ and
$v_1\mcom v_n$ which satisfies $\epsilon(\tau)=\s{v_k v_{k+1}}$
when $\tau=(k,k+1)$ is a transposition.

In addition to the algebra structure, $\bigwedge V$ also possesses
a natural \ztz-cocommutative coalgebra structure given by
\begin{multline*}
\Delta(v_1\mwedge v_n)=\\
\sum_{k=1}^n \sum_{\sigma\in \sh(k,n-k)}
\s{\sigma}\epsilon(\sigma) v_{\sigma(1)}\mwedge v_{\sigma(k)}\tns
v_{\sigma(k+1)}\mwedge v_{\sigma(n)},
\end{multline*}
where $\sh(k,n-k)$ is the set of all unshuffles of type $(k,n-k)$,
that is, permutations which are increasing on $1\dots k$, and on
$k+1\dots n$.

Note that in our definition of the exterior algebra, we begin with
elements of degree 1, so that  $\k=\bigwedge^0 V$ is not
considered part of the exterior algebra. Therefore, in our
construction $\Delta$ is not injective, because $\Delta(v)=0$ for
any $v\in V$; in fact, $\ker \Delta= V$.

In order to generalize Lie algebras to \linf\ algebras, it is
convenient to replace the space $V$ with its parity reversion
$W=\Pi V$, and the exterior coalgebra on $V$ with the symmetric
coalgebra on $W$. One advantage of this construction is that the
symmetric coalgebra $\SW$ is \zt-graded, rather than \ztz-graded,
and the coderivations on $\SW$ determine a \zt-graded Lie
algebra, which is better behaved in terms of the properties of
the bracket than the \ztz-graded coderivations on $\bigwedge V$.
In case the space $V$ is not a \zt-graded space, that is to say
that all elements in $V$ have even parity, then $W$ becomes a
space consisting of only odd elements, and the symmetric algebra
on $W$ coincides in a straightforward manner with the exterior
algebra on $V$, in the sense that the exterior degree of an
element in $\bigwedge V$ is the same as the parity of its image
in $\SW$. Thus in the classical picture one studies simply the
space $\bigwedge V$ equipped with the induced \zt-grading given
by the \Z-grading inherited from the tensor algebra.

Let us denote the product in $\SW$ by juxtaposition. Then the
coalgebra structure on $\SW$ is given by the rule
\begin{equation*}
\Delta(w_1\cdots w_n)= \sum_{k=1}^n \sum_{\sigma\in \sh(k,n-k)}
\epsilon(\sigma) w_{\sigma(1)}\cdots w_{\sigma(k)}\tns
w_{\sigma(k+1)}\cdots w_{\sigma(n)}.
\end{equation*}
Again the kernel of $\Delta$ is simply $W$, and it is injective on
elements of higher degree.

There is a natural correspondence between $\hom(\SW,W)$ and the
\zt- graded Lie algebra of coderivations $\coder(W)$, which is
given by extending a map $\ph:S^k(W)\ra W$ to a \zt-graded
coderivation by the rule
\begin{equation*}
\ph(w_1\cdots w_n)=\sum_{\sigma\in\sh(k,n-k)}\epsilon(\sigma)
\ph(w_{\sigma(1)}\cdots w_{\sigma(k)})w_{\sigma(k+1)}\cdots
w_{\sigma(n)}.
\end{equation*}
Moreover, if $\ph$ is any coderivation, and $\ph_k:S^k(W)\ra W$
denotes the induced maps, then $\ph$ can be recovered from these
maps by the formula
\begin{equation*}
\ph(w_1\cdots w_n)=\sum_{1\le k\le n}
\sum_{\sigma\in\sh(k,n-k)}\epsilon(\sigma)
\ph_k(w_{\sigma(1)}\cdots w_{\sigma(k)}) w_{\sigma(k+1)}\cdots
w_{\sigma(n)}.
\end{equation*}
We thus can express $\ph=\sum_{n=1}^\infty \ph_n$, and say that
$\ph_n$ is the degree $n$ part of $\ph$. Note that $\coder(W)$ is
actually a direct product, rather than direct sum, of its graded
subspaces.  In the case of Lie and associative algebras, it is
conventional to consider the direct sum of the subspaces instead,
but we shall not make that convention here.

If $\ph$ and $\psi$ are two coderivations, then their Lie bracket
is given by
\begin{equation*}
[\ph,\psi]=\ph\circ\psi-\s{\ph\psi}\psi\circ\ph.
\end{equation*}
In terms of the identification of $\coder(W)$ with $\hom(\SW,W)$,
this bracket takes the form
\begin{multline*}
[\ph,\psi]_n(w_1\cdots w_n)=\\
\sum_{k+l=n+1\atop\sigma\in\sh(k,n-k)} \epsilon(\sigma)
[\ph_l(\psi_k(w_{\sigma(1)}\cdots w_{\sigma(k)})
w_{\sigma(k+1)}\cdots w_{\sigma(n)})\\
-\s{\ph\psi} \psi_l(\ph_k(w_{\sigma(1)}\cdots w_{\sigma(k)})
w_{\sigma(k+1)}\cdots w_{\sigma(n)}].
\end{multline*}

A codifferential $\d$ is a coderivation whose square is zero, and
in case $\d$ is odd, this is equivalent to the property
$\br\d\d=0$, which can be expressed in the form
\begin{equation*}
\sum_{k+l=n+1\atop\sigma\in\sh(k,n-k)} \epsilon(\sigma)
[d_l(d_k(w_{\sigma(1)}\cdots w_{\sigma(k)}) w_{\sigma(k+1)}\cdots
w_{\sigma(n)})]=0.
\end{equation*}

Let $\L_n=\hom(S^n(W),W)$, and $L=\prod_{n=1}^\infty
L_n=\coder(W)$. An odd codifferential $d_2$ in $\L_1$, called a
quadratic codifferential, determines a \zt-graded Lie algebra
structure on $V$. The reader can check that the codifferential
property is precisely the requirement that the induced map
$l_2:V\wedge V\ra V$ satisfies the \zt-graded Jacobi identity. In
order to define a \linf\ algebra, one simply takes an arbitrary
codifferential in $L$. This codifferential $d$ is given by a
series of maps $d_k:S^{k}(W)\ra W$, which determine maps
$l_k:\bigwedge^{k}V\ra V$, satisfying some relations which are
the defining relations of the \linf\ algebra. Details of this
construction can be obtained in \cite{pen1}.

It is easy to see that the map $D:\L\ra\L$ given by
$D(\ph)=[d,\ph]$ is a differential, equipping $L$ with the
structure of a differential graded Lie algebra. When $d$ is
quadratic, the homology given by this differential is essentially
the Eilenberg-Chevalley cohomology of the Lie algebra (with
coefficients in the adjoint representation). In general, the
homology given by $d$ is called the cohomology of the \linf\
algebra $V$. Thus both Lie and \linf\ algebras are structures on a
\zt-graded vector space, which are determined as odd
codifferentials on the symmetric coalgebra of the parity
reversion of the space, and the cohomology of these algebras is
simply the homology on the graded Lie algebra of coderivations of
this coalgebra determined by this codifferential.

A similar picture applies for associative algebras, except that
this time, we form the tensor coalgebra $\TW$ of the parity
reversion $W$ of the underlying vector space $V$, for which we
take the \zt-grading corresponding to the grading of $W$, instead
of the \ztz-grading which would be necessary if we were to
consider the coalgebra structure on $\TV$. The coproduct is given
by
\begin{equation*}
\Delta(w_1\mcom w_n)=\sum_{k=1}^{n-1}(w_1\mcom w_k)\tns
(w_{k+1}\mcom w_n).
\end{equation*}

There is a natural correspondence between $\hom(\TW,W)$ and the
\zt- graded Lie algebra of coderivations $\coder(W)$, which is
given by extending a map $\ph:T^k(W)\ra W$ to a \zt-graded
coderivation by the rule
\begin{multline*}
\ph(w_1\mcom w_n)=\sum_{0\le i\le n-k}[\\
\s{(w_1\mplus w_i)\ph} (w_1\mcom w_i,\ph(w_{i+1}\mcom w_{i+k}),
w_{i+k+1}\mcom w_n)].
\end{multline*}
Moreover, if $\ph$ is any coderivation, and $\ph_k:T^k(W)\ra W$
denotes the induced maps, then $\ph$ can be recovered from these
maps by the formula
\begin{multline*}
\ph(w_1\mcom w_n)=\\
\sum_{1\le k\le n\atop 0\le i\le n-k} \s{(w_1\mplus w_i)\ph}
(w_1\mcom w_i,\ph_k(w_{i+1}\mcom w_{i+k}), w_{i+k+1}\mcom w_n) .
\end{multline*}
If $\ph$ and $\psi$ are two coderivations, then their Lie bracket,
in terms of the identification of $\coder(W)$ with $\hom(\TW,W)$,
takes the form
\begin{multline*}
[\ph,\psi]_n(w_1\mcom w_n)=
\sum_{k+l=n+1\atop 0\le i\le n-k}[\\
\s{(w_1\mplus w_i)\psi}
\ph_l(w_1\mcom w_i,\psi_k(w_{i+1}\mcom w_{i+k}), w_{i+k+1}\mcom w_n)-\\
\s{(\psi+w_1\mplus w_i)\ph} \psi_l(w_1\mcom
w_i,\ph_k(w_{i+1}\mcom w_{i+k}), w_{i+k+1}\mcom w_n)].
\end{multline*}

An odd codifferential $d$ on $\TW$ satisfies the relations
\begin{multline*}
\sum_{k+l=n+1\atop 0\le i\le n-k} \s{w_1\mplus w_i}d_l(w_1\mcom
w_i,d_k(w_{i+1}\mcom w_{i+k}),w_{i+k+1}\mcom w_n)\\=0.
\end{multline*}
The relations satisfied by the induced maps $m_k:V^k\ra V$
determine the structure of an \ainf\ algebra on $V$. When $d$ is
a quadratic codifferential, one can show that the induced map $m$
is simply an associative algebra structure on $V$. As in the case
of an \linf\ algebra, one forms the graded Lie algebra
$L=\prod_{n=1}^\infty L_n$, where $L_n=\hom(T^{n}(V),V)$ is the
subspace consisting of degree $n$ elements in $L$. The
differential $D(\ph)=[\ph, d]$ equips $L$ with the structure of a
differential graded Lie algebra, and the homology given by $D$ is
called the cohomology of the \ainf\ algebra.  For an associative
algebra, the cohomology is simply the Hochschild cohomology of
the algebra, and the bracket is the Gerstenhaber bracket.

At this point, it should be clear to the reader that from an
abstract point of view there is not much difference between the
construction of \linf\ and \ainf\ algebras. In both cases we have
a graded Lie algebra $L$ consisting of \zt-graded coderivations
of an appropriate coalgebra, which can be expressed as a direct
product of subspaces $\L_n$, for $n\ge 0$.  It is easy to see
that $[\L_n,\L_m]\subseteq \L_{m+n-1}$. From this it follows that
if $d$ is a quadratic codifferential, then $D(L_n)\subseteq
L_{n+1}$, so that the homology $\HL$ has a decomposition
$\HL=\prod_{n=1}^\infty H_n(\L)$, where $H_n(\L)=\ker(d:L_n\ra
L_{n+1})/\im(d:\L_{n-1}\ra \L_n)$. When $d$ is not quadratic, no
such decomposition exists; instead we need to consider a natural
filtration on $\HL$. The filtration arises from the natural
filtration on $\L$ given by $L^n=\prod_{i=n}^\infty L_i$, which
is respected by $D$ in the sense that $D(\L^n)\subseteq L^n$. We
shall discuss how the filtration on $\L$ gives rise to a natural
filtration on $\HL$ in the next section.

For Lie and associative algebras, infinitesimal deformations are
determined by the second cohomology group $H^2(\L)$. For infinity
algebras, the entire cohomology governs the infinitesimal deformations.
In fact, one interpretation of the cohomology of a Lie or associative
algebra is that it governs the infinitesimal deformations of the
algebra into the appropriate infinity algebra.

Suppose now that $\A$ is a \zt-graded commutative algebra over a
ground field \k\ equipped with a fixed augmentation
$\epsilon:\A\ra\k$, with augmentation ideal $\m=\ker(\epsilon)$.
For simplicity, denote $\WA=\W\tns\A$. Let $\T_{\A}(\WA)$ be the
tensor algebra of $\WA$ over \A. Then
$\T_{\A}(\WA)\cong\TW\tns\A$. Similarly, if $S_{\A}(\WA)$ is the
symmetric algebra of $\WA$ over $\A$, then
$S_{\A}(\WA)\cong\SW\tns\A$. These natural isomorphisms respect
the algebra and coalgebra structures of both sides. In what
follows, let us assume for sake of definiteness that we are
working with the symmetric coalgebra, and thus with the
deformation theory for \linf\ algebras, but the statements and
results hold true for \ainf\ algebras as well, by simply replacing
the symmetric coalgebra with the tensor coalgebra.

Let $\LA=\L\htns\A=\prod_{n=1}^\infty L_n\tns A$ be the completed
tensor product of $L$ and \A. There is a natural identification
of \LA\ with the space of \A-algebra coderivations of $\SWA$.
Both $\LA$ and $\L$ have the structure of \zt-graded \A-Lie
algebras, where the augmentation determines the \A-Lie algebra
structure on $L$. The projection $\LA\ra\L$ induced by the
augmentation is an \A-Lie algebra homomorphism.

A \emph{deformation of a \linf\ structure $d$ with base \A}, or
more simply an \A-deformation of $d$, is defined to be an odd
codifferential $\td\in\LA$, which maps to $\d$ under the natural
projection.  When $\d\in\L_2$, so that the structure is a Lie
algebra, then an \A-deformation of the Lie algebra structure is
given by an odd codifferential in $(\LA)_2$.

The algebra \A\ splits canonically in the form $\A=\k\oplus\m$,
which induces a splitting
$\LA=\L\tns\k\oplus\L\htns\m=L\oplus\Lm$; moreover, the map
$\LA\ra\L$ is simply the projection on the first factor. Thus if
$\td$ is an \A-deformation of $L$ then $\td=\d+\delta$ where
$\delta\in\Lm$. In order for $\td$ to be an \A-deformation, we
must have $[\td,\td]=0$, which is equivalent to $\delta$
satisfying the \emph{Maurer-Cartan formula}
\begin{equation}\label{mc formula}
D(\delta)=-\tfrac12\br\delta\delta,
\end{equation}
where $D(\delta)=\br\d\delta$.

 If $\lambda:W\ra W'$ is even, then it extends uniquely to a coalgebra
homomorphism $S(\lambda):\SW\ra S(W')$ by
\begin{equation*}
S(\lambda)(w_1\cdots w_n)=\lambda(w_1)\cdots\lambda(w_n).
\end{equation*}
If $W$ and $W'$ are equipped with Lie algebra structures $\d$ and
$d'$ respectively, then  $\lambda$ determines a Lie algebra homomorphism
if $\d'\circ S(\lambda)=S(\lambda)\circ\d$.  For the \linf\ case,
a homomorphism is given by an arbitrary coalgebra morphism
$f:\SW\ra S(W')$ satisfying $\d'\circ f=f\circ d$, where $d$ and
$d'$ are now \linf\ algebra structures. It is not possible to
translate the definition of homomorphism of \linf\ or Lie
algebras into a statement about the algebras $\L(W)$ and
$\L(W')$, but this is not surprising since homomorphisms of Lie
algebras do not induce morphisms on the cohomology level.

If $f:\SW\ra \SW$ is a coalgebra automorphism, then $f$ induces
an automorphism $f^*$ of $\L$, given by
$f^*(\ph)=f\inv\circ\ph\circ f$, so that $d_f=f^*(\d)$ is a
codifferential in $L$ if $d$ is a codifferential. If $f=S(\tau)$
for some isomorphism $\tau:\W\ra\W$,  then $f^*(\L_n)=\L_n$.
In particular, $d_f$ is a quadratic codifferential when $d$
is. For Lie algebras, we restrict our consideration of
automorphisms to these degree zero maps, but for \linf\ algebras
we require merely that $f$ be even, so that it may mix the
exterior degrees of elements.

Two \A-deformations $\td$ and $\td'$ of an \linf\ algebra are
said to be \emph{equivalent} when there is an \A-coalgebra
automorphism $f$ of $\SWA$ such that $f^*(\td)=\td'$, compatible
with the projection $\SWA\ra\SW$. For the Lie algebra case, one
requires in addition that $f=S(\gamma)$ for some isomorphism
$\gamma:\WA\ra\WA$. In terms of the decomposition
$\SWA=\SW\bigoplus\SW\tns\m$ we can express an equivalence in the
form $f=\id+\lambda$, for some map $\lambda:\SW\ra\SW\tns\m$.
Since $f$ must satisfy the condition
\begin{equation*}
\Delta\circ f= f\tns f\circ\Delta,
\end{equation*}
and $\id$ is an automorphism, $\lambda$ must satisfy the condition
\begin{equation}
\label{autocond} \Delta\circ\lambda=(\lambda\tns\id
+\id\tns\lambda+\lambda\tns\lambda)\circ\Delta.
\end{equation}
We can also express $f^*=\id+\tl$, for some map $\tl:\L\ra\Lm$,
and we obtain the condition
\begin{equation}
\label{Lautocond} \tl\br{\vphantom{\tl\ph}\ph}\psi=
\br{\tl\ph}{\psi}+\br{\ph}{\tl\psi}+\br{\tl\ph}{\tl\psi}.
\end{equation}

\vspace{10pt}

We next define the notion of a change of base of a deformation.
Suppose that $\A$ and $\A'$ are two augmented \k-algebras with
augmentation ideals $\m$ and $\m'$ resp. and
 $\tau:\A\ra\A'$ is a \k-algebra morphism.
Then $\tau$ induces an \A-linear map
$\tau_*=1\tns\tau:\SWA\ra\SWA'$ which is an \A-coalgebra
morphism, that is
\begin{equation*}
\Delta'\circ\tau_*=(\tau_*\tns\tau_*)\circ\Delta.
\end{equation*}
Similarly, the induced map $\tau_*:\LA\ra\LA'$ is a homomorphism
of graded \A-Lie algebras. In terms of the decompositions
$\LA=\L\bigoplus\Lm$ and $\LA'=\L\bigoplus\Lm'$, it is clear that
$\tau_*(\ph)=\ph$ for any $\ph\in\L$, and
$\tau_*(\Lm)\subseteq\Lm'$. If $\td=\d+\delta$ is an
\A-deformation of $\d$, then the {\em push out}
$\tau_*(\td)=\d+\tau_*(\delta)$ is an $\A'$-deformation of $\d$.
Furthermore, if two \A-deformations of $\d$ are equivalent, then
$\tau_*$ maps them to equivalent deformations, so that $\tau_*$
induces a map between equivalence classes of \A-deformations  and
equivalence classes of $\A'$-deformations.

Two special cases of deformations with base arise, the
infinitesimal deformations and the formal deformations. If \A\ is
a local algebra and $\m^2=0$ then we shall call \A\ an
infinitesimal algebra, and an \A-deformation will be called an
infinitesimal deformation. Especially interesting is the
classical notion of an infinitesimal deformation, determined by
the algebra $\A=\k[t]/(t^2)$, where $t$ is taken as an even
parameter, which may be generalized to the dual algebra
$\A=\k[t,\theta]/(t^2,t\theta,\theta^2)$, where $\theta$ is taken
as an odd parameter.

A formal deformation is given by taking $\A$ to be a complete
local algebra, or a formal algebra, so that
$\A=\invlim_{n\ra\infty} A/\m^n$. Then a \emph{formal}
deformation with base \A\ is given by a codifferential $\td$ on
$\invlim_{n\ra\infty}\L\hat\tns\A/\m^n$, the classical example
being $\A=\k[[t]]$,  the ring of formal power series in the even
parameter $t$.

The main purpose of this section was to show how the notion of an
\A-deformation  can be given in terms of the Lie algebra $\LA$.
Deformations with base \A\ are given by certain classes of
codifferentials in \LA, depending on the type of structure being
deformed. Equivalences of deformations are given by certain classes of
automorphisms of the Lie algebra structure of \LA, again depending on
the type of structure being deformed. We will show later how the
homology of $L$ determined by the codifferential $\d$ relates the
notion of infinitesimal deformation with that of infinitesimal
equivalence.

\section{Filtered Topological Vector Spaces}
Because we will need some topological results, we include a brief
explanation of properties of \emph{filtered topological modules}.
A filtered topological space is a natural generalization of a
direct product, and arises in our construction because the space $L$
of cocycles does not have a natural decomposition as a sub-direct
product of the space of cochains. The natural direct product
decomposition of $L$ does induce a filtration on the space of cocycles,
which descends to a filtration on the homology.  The topology of 
a filtered space allows the natural introduction of a dual space,
the continuous dual, which is small enough to be useful in our
construction.

A (decreasing) filtration on a module $F$ is a sequence of submodules $F^n$
satisfying $F=F^0$, $F^{n+1}\subseteq F^n$ and
$\bigcap_{n=1}^\infty F^n=0$.  A sequence $\{x_i\}$ is said to be
Cauchy if given any $n$, there is some $m$ such that if $i,j\ge
m$, then $x_i-x_j\in F^n$. The space $F$ is said to be complete
if every Cauchy sequence converges, that is if there is some $x$
such that for any $n$ there is some $m$ such that $x-x_i\in F^n$
if $i\ge m$. This $x$ is unique, and is called the limit of the
sequence $\{x_i\}$. This notion of convergence, whether or not
$F$ is complete, determines a topology on $F$, which we will call the
filtered topology. The order of an element $x$ is the largest $k$
such that $x\in F^k$. The only element which has infinite order
is zero.

A map $f:F\ra G$ of two filtered topological spaces is continuous
iff for given $n$ there is some $m$ such that $f(x)\in G^n$ whenever
$x\in F^m$. To see this, note that $f$ is continuous iff
$\{f(x_i)\}$ is a Cauchy sequence whenever $\{x_i\}$ is Cauchy.

Let $F_i$ be a subspace of $F^i$ which projects isomorphically to
$F^i/F^{i+1}$. Then $\bar F=\prod_{i=0}^\infty F_i$ has a natural
filtration ${\bar F}^n=\prod_{i=n}^\infty F_i$, and is a complete
filtered topological space. There is a natural map
$\iota:F\ra\bar F$ defined as follows. Let $x\in F$. Then there
is a unique $x_0\in F_0$ whose image in $F_0/F_1$ coincides with
that of $x$. Then $x-x_0\in F^1$. Continuing, one obtains a
sequence of elements $x_i\in F_i$ such that $x-\sum_{i=0}^n
x_i\in F^{n+1}$. Define $\iota(x)=\prod x_i\in\bar F$. The
natural map $\iota$ is injective, order preserving, and
continuous, and is surjective precisely when $F$ is complete. In
this case, the inverse map $\iota\inv$ is also continuous, so that
$F\cong\bar F$. From this we see that a complete filtered
topological space is essentially the same as a direct product.

If $F=\prod F_i$ is complete, and $\{x^n_i\}$ is a subset of
$F^n$ which projects to a basis of $F^n/F^{n+1}$, then any
element $x$ in $F$ has a unique expression as an infinite sum
$x=a_n^ix^n_i$ where $a_n^i\in\k$ and for fixed $n$, only
finitely many of the coefficients $a^i_n$ are non zero. Note that
$x^n_i$ has order $n$. An ordered set $\{y_i\}$ which satisfies
the property that every element of $F$ can be written uniquely as
an infinite sum $a^iy_i$ will be called a basis of $F$, it is
increasing if $o(y_n)\ge o(y_{n+1})$, and strictly increasing if
for any $n$ there is some $m$ such that $\{y_i\}_{i\ge m}$ is a
basis of $F^n$. The basis $x^n_i$ can be ordered in a strictly
increasing manner.

If $B$ is a subspace of $F$ then it inherits a natural filtration
$B^i=B\cap F^i$. If $F$ is complete then the subspace $B$ is
closed in $F$ precisely when it is complete as a filtered space.
The space $H=F/B$ can be given a filtration by $H^i=\rho(F^i)$
where $\rho:F\ra F/B$ is the canonical map. One obtains
$H^i=F^i/B^i$ in a natural way. The requirement
$\bigcap_{i=0}^\infty H^i=0$ is satisfied precisely when $B$ is
closed in $F$. For suppose that $B$ is not closed, and $\{b_i\}$
is a sequence in $B$ converging to some $b\not\in B$. Let
$x_i=b-b_i$, so that $\rho(x_i)=\rho(b)$ for all $i$, but this
implies that $\rho(b)\in H^i$ for all $i$. On the other hand,
suppose that $h\ne 0\in H^i$ for all $i$. Then $h=f(x_i)$ for
some $x_i\in F^i$. Let $y_i=x_0-x_{i}$. Then $\rho(y_i)=0$ for
all $i$, so $y_i\in B$, but the sequence $y_i$ converges to
$x_0$, which does not lie in $B$. Thus the quotient space $F/B$
is a filtered space precisely when $B$ is closed in $F$.
Suppose that $F$ is complete and $B$ is closed in $F$. Then we
claim that $H=F/B$ is also complete. For suppose that $h_i$ is a
Cauchy sequence in $H$, and by taking a subsequence if necessary,
we may assume that $h_i-h_{i+1}\in H^i$. Choose $x_i\in F^i$ such
that $\rho(x_i)=h_i-h_{i+1}$, and $a\in F$ such that
$\rho(a)=h_0$. Let $y_n=a-\sum_{i=0}^n x_i$. Then
$y_n-y_{n+1}=x_{n+1}\in F^n$, so $y_n$ converges to some $y$.
Note that $\rho(y_n)=h_{n+1}$, so it follows that $h_i$ converges
to $\rho(y)$.

Next, suppose that $f:F\ra G$ is a continuous map of filtered
spaces, $Z=\ker f$ and $B=\im f$. Then both $Z$ and $B$ inherit
the structure of filtered spaces. It is easy to see that $Z$ is
closed in $F$, but it may happen that $B$ is not closed in $G$.
For example, if $F$ is a filtered space which is not complete,
then its image in its completion under the canonical injection is
not closed. Completeness of $F$ is also not sufficient to
guarantee that the image is closed, as can be seen from the
following example.  Let $F_0$ have a countably infinite basis
$x_i$, and $F_1=0$.  Then $F$ is complete for trivial reasons.
Let $y_i$ be a basis for $G_i$ and $G=\prod G_i$ with the natural
filtration. Then $G$ is complete. Define a map $f:F\ra G$ by
$f(x_i)=y_i$. Then $f$ is continuous, but its image is not closed
in $G$. Let us say that a filtered space is of finite type if
$\dim(F^n/F^{n+1})<\infty$ for all $n$.

\begin{lma}
If $F$ is a complete filtered space of finite type, and $f:F\ra
G$ is continuous, then $B=\im f$ is closed in $G$.
\end{lma}
\begin{proof}
We may as well assume that $F=\prod F_i$ with the standard
filtration, where each $F_i$ is finite dimensional. If we take a
basis $\{x_i^n\}$ of $F_n$, we can obtain a strictly increasing
basis of $F$, and by throwing out unnecessary elements, we can
choose a subsequence which spans a subspace mapping injectively
to $B$. Let $\{y_i\}$ be the subset of $B$ so obtained, ordered
in a strictly increasing manner, because by continuity, given any
$n$, there can only be a finite number of the $x_i^m$ whose image
has order smaller than $n$. Next we claim that an element $y$
lies in $B$ precisely if it can be written as an infinite sum of
the form $y=a^ky_k$. To see this fact, first note that because
the order of the $x_i^n$'s are increasing, and $y_k$ is the image
of some $x_{i_k}^{n_k}$, we can form the element
$a^kx_{i_k}^{n_k}$ which is well defined in $F$ since it is
complete, and the image of this element must be $y$. On the other
hand, by our construction, if $y=f(a_n^ix_i^n)$, then by using
the fact that for any finite combination of the $f(x_i^n)$ can be
expressed in terms of the $y_i$, we can subtract off a linear
combination $b^i y_i$ from $y$ such that $y-b^iy_i$ is expressed
as an image of terms of high order in $F$. This therefore must
have high order in $G$, and thus we can express $y$ as the limit
of a Cauchy sequence in the $y_i$.
\end{proof}

Now let us suppose that $F$ has a coboundary operator $\d$ which
respects the filtration, that is to say that $d^2=0$ and
$d(F^n)\subseteq F^n$ for all $n$. Then if $F$ is of finite type,
the homology $H(F)$ has the natural structure of a filtered
space, and is complete if $F$ is complete. In the case where $F$
is a direct product of finite dimensional vector spaces, one
immediately obtains that the homology is a complete filtered
space.

Consider $\k$ as a filtered space with $\k^1=0$. Then we can form
the continuous dual space $F^*$. It consists of all continuous
linear functionals, that is all $\lambda:F\ra\k$ such that there
is some $n$ such that $\lambda(x)=0$ for all $x\in F^n$. 
The (continuous) dual of a filtered space is not filtered in the sense we
have described above, but does possess an increasing filtration.
To distinguish this type of space from the filtered spaces 
given by decreasing filtrations, let us say
that a space $E$ is \emph{\gr} if there is a sequence $E_i$ of subspaces
satisfying $E_0=\{0\}$, $E_n\subseteq E_{n+1}$, and $\bigcup E_n=E$. 
If $F$ is filtered, then $F^*$ is \gr, with
$(F^*)_n=\{\lambda\in F^*|\lambda(F^n)=0\}$.  As a filtered space is a model
of a direct product, a \gr\ space is a model of a direct sum.  For
suppose that we choose subspaces $E^k\subseteq E_{k+1}$ such that
$E^k$ projects isomorphically to $E_{k+1}/E_k$. Then there is a
natural isomorphism $E\ra\oplus E^k$.
(This accounts for our terminology. The usual definition of a graded
space is one given by a decomposition as a direct sum as above. Our
definition generalizes the usual one.)
 We do not equip a \gr\ space
with any topology, but still, its dual space has a natural filtration
given by $(E^*)^n=\{\lambda\in E^*|\lambda(E_n)=0\}$. Moreover, the
dual of a \gr\ space is complete.  

Let us say that a
$f:E\ra D$ map between two \gr\ spaces respects the grading if given
$n$, there is some $m$ such that $f(E_n)\subseteq E_m$. 
A graded space is said to be of \emph{finite type}
if all the $E_n$ are finite dimensional, in which case every
map from $E$ to a graded space respects the grading. 
If $f:F\ra G$ is a continuous map of filtered spaces, then it induces
a map $f^*:G^*\ra F^*$, which respects the grading,
giving a contravariant functor from the category
of filtered spaces with continuous maps to the category of \gr\ spaces
with maps respecting the grading. If $f:E\ra D$ respects the grading,
then $f^*:D^*\ra E^*$ is continuous, so we again a contravariant
functor between graded spaces and filtered spaces.

It is useful  to note that if $D$ is a subspace of a graded space, then
$D$ inherits a grading given by $D_i=D\cap E_i$, $E/D$ is 
graded by $(E/D)_i=\im(E_i)$, and the inclusion and projection maps
respect the grading.  Furthermore, any subspace of a graded space has
a graded complementary subspace.

For filtered and \gr\ spaces, a space is of finite type precisely when
its dual is of finite type. Moreover, if $F$ is a finite type filtered space, 
then $(F^*)^*$ is naturally identified with $\bar F$.  For \gr spaces,
the tensor product $E\tns D$ can be graded by $(E\tns
D)_n=\sum_{p+q= n}E_p\tns E_q$ (the sum is not direct). Similarly, the
tensor product $F\tns G$ of filtered spaces has a filtration
$(F\tns G)^n=\sum_{p+q=n} F^p\tns G^q$. Then one obtains the useful 
formulae $(E\tns D)^*=E^*\htns F^*$ and $(F\tns G)^*=F^*\tns G^*$. 

When $F$ is a finite type filtered space, and
$\{x_i\}$ is a strictly increasing basis of $F$, then the dual
basis $\{\lambda^i\}$, given by $\lambda^i(x_j)=\delta_j^i$, is a
well defined basis of $F^*$.  It is clear that $\lambda^i$ is
continuous, so we only need to show that any continuous linear
functional $\lambda$ can be represented as a sum of the
$\lambda_i$. Let $a_i=\lambda(x_i)$. By continuity, the sum
$a_i\lambda^i$ is well defined and coincides with $\lambda$.

Now suppose that $M$, $N$ are filtered spaces, $N$ has an
increasing basis $\{y_i\}$, while $M$ is of finite type
and has an increasing basis
$\{x_i\}$ with dual basis $\lambda^i$ of $M^*$.  Let $\eta$ be an
element of $N\tns M^*$ which can be written in the form $\eta=a^i_j
y_i\tns \lambda^j$ for some finite sequence of elements
$a^i_j\in\k$. Clearly $\eta$ determines a continuous map $M\ra N$
by the rule $\eta(b^kx_k)=a^i_kb^ky_i$. Introduce the filtration
on $N\tns M^*$ induced by the filtration on $N$.  Then one can
form the completion $N\htns M^*$ of this filtration, which will
have basis $\{y_i\tns\lambda^j\}$. Any element $\eta$ of the
completion has a unique expression in the form $\eta=y_i\tns
\beta^i$, where $\beta^i\in M^*$ (of course, only a finite number
of terms of each order can occur). When $N$ is complete and of
finite type, then $\eta$ also determines a continuous map from
$M$ to $N$, and moreover, any continuous map is so obtained.
Thus we can identify $N\tns M^*$ with $\hom(M,N)$.  In general, this
filtration will not be of finite type even when $N$ has finite
type and is complete. Our main interest will be in representing
elements of $\hom(M,N)$ by elements in $N\htns M^*$. When $M$ is
of finite type, and $N$ is not necessarily complete, we still have
\begin{equation*}
\hom(M,N)\subseteq N\htns M^*\subseteq\bar N\htns M^*=\hom(M,\bar
N).
\end{equation*}
In particular, every continuous linear map has a representation
as an element of $N\htns M^*$. Finally, note that when $M$ is of
finite type, and $m^i$ is an increasing basis of $M^*$, then any
element of $N\htns M^*$ can be expressed uniquely in the form
$n_i\tns m^i$, where $n_i$ is an increasing sequence in $N$.

\section{Infinitesimal Deformations}

Let $\A$ be a \zt-graded commutative algebra, $L$ a differential graded Lie
algebra and
$\td=\d+\delta$ be an infinitesimal \A-deformation of $\L$.
Since $\m^2=0$, the Maurer-Cartan formula (equation \ref{mc
formula}) reduces to the cocycle condition $D(\delta)=0$.
Furthermore, an infinitesimal equivalence is of the form
$f=1+\lambda$ where $\lambda$ is a coderivation of $\SWA$,
because equation \ref{autocond} reduces to
$\Delta\circ\lambda=(\lambda\tns\id
+\id\tns\lambda)\circ\Delta$.  Furthermore, if we express
$f^*=\id+\tl$, then $\tl\ph=\br\ph\lambda$, since
$f\inv=\id-\lambda$, so that \begin{equation*}
f^*(\ph)=(\id-\lambda)\circ\ph\circ(\id+\lambda)= \ph
+\ph\circ\lambda-\lambda\circ\ph=\ph+\br\ph\lambda.
\end{equation*}
A trivial infinitesimal deformation is one of the form
$f^*(d)=d+D(\lambda)$. Thus the (even part of the) homology
$\HL\tns\m$ classifies the equivalence classes of infinitesimal
\A-deformations of an \linf\ algebra. In the case of a Lie
algebra, the derivations $\lambda$ giving rise to infinitesimal
automorphisms are all determined by linear maps, so are elements
of $\L_0\tns\m$, while the infinitesimal deformations are given by
elements of $\L_1\tns\m$, so that $H_1(\L)\tns\m$ classifies the
equivalence classes of deformations. If $\d+\delta$ and
$\d+\delta'$ are two \A-deformations of $\d$, then they are
equivalent precisely if $\delta-\delta'$ is a coboundary, while
the condition for $\d+\delta$ to be an \A-deformation is simply
that $\delta\in\ZL\htns\m$, where $\ZL=\ker D$ is the space of
cocycles.

If $\tau:\A\ra\A'$ is a morphism of infinitesimal \k-algebras,
then $\tau_*:\LA\ra\LA'$ induces a morphism of equivalence classes
of infinitesimal deformations of some fixed codifferential
$\d\in\L$. A universal infinitesimal deformation of $\d$ is an
initial object in the category of such equivalence classes. In
other words, a universal infinitesimal deformation of $\d$ is
given by some infinitesimal \k-algebra \A\ and \A-deformation
$\td=\d+\delta$, such that if $\A'$ is another infinitesimal
\k-algebra, and $\td'=\d+\delta'$ is an $\A'$-deformation of
$\d$, then there is a unique morphism $\tau$ of infinitesimal
\k-algebras satisfying the property that $\tau_*(\td)$ is
equivalent to $\td'$.

For a \zt-graded space, the parity reversion of the space is
given by interchanging the parity of homogeneous elements, so
even elements give rise to odd elements in the parity reversion
and vice versa. Let $\Pi \HL$ denote the parity reversion of
$\HL$, and equip it with the filtration which it inherits from
$\HL$. Let $\m=(\Pi\HL)^*$. Let $\A=\k\bigoplus\m$, with
multiplication on $\m$ defined trivially, so that $\A$ is an
infinitesimal \k-algebra. Let $\mu:\Pi\HL\ra\ZL$ be a right
inverse to the map $\ZL\ra\Pi\HL$, respecting the filtrations on
$\Pi\HL$ and $\ZL$, in other words, $\mu(\pi\delta)\in\ZL^n$ if
$\delta\in\HL^n$. We want to represent $\mu$ as an element of
$\ZL\hat\tns\m$.

Let us assume that $\HL$ is of finite type, and choose a strictly
increasing basis $\delta_i$ of $\HL$, Define $t^i\in(\Pi\HL)^*$
to be the dual basis, \ie, $t^i(\pi\delta_j)=\delta_j^i$. If we
let $\mu_i=\mu(\pi\delta_i)$, then $\mu_i\htns t^i$ represents
the map $\mu$.

The reason we work with the parity reversion $\Pi\HL$ instead of
$\HL$ is because we want $\mu$ to be an odd map, because we are
going to define an \A-deformation $\d+\mu=\d+\mu_i\htns t^i$. So
the parity reversion does the trick, because it turns an even map
$\HL\ra\ZL$ into an odd map $\Pi\HL\ra\ZL$. Also it is immediate
that $\d+\mu$ is an \A-deformation, because $D(\mu)=0$, since
$\mu_i\in\ZL$ for all $i$.

Moreover, if $\mu'$ is another inverse map of the canonical map
$\ZL\ra\HL$, and we express $\mu'=\mu'_i\tns t^i$, then
$\mu_i-\mu_i'=\delta(\ph_i)$ for some $\ph_i$. Thus we obtain that

\begin{equation*}
\mu-\mu'=(\mu_i-\mu_i')\tns t^i =D(\ph_i)\tns t^i=D(\ph_i\tns
t^i)
\end{equation*}

is a coboundary, so that the two \A-deformations are
equivalent. Note that since the order of $\mu_i-\mu_i'$ is
increasing, we can assume that the sequence $\{\ph_i\}$ has only
finitely many terms of any order, and therefore $\ph_i\tns t^i$
is a well defined element of $\L\tns\m$.

Now suppose that $\A'$ is another infinitesimal algebra with
augmentation ideal $\m'$, and $\d+\delta'$ is an arbitrary
$\A'$-deformation. Then $\delta$ can be expressed in the form
$\delta'=\mu_i\tns m^i+ b$, where $ b\in\BL$, and by replacing it
with an equivalent deformation we can assume that $b=0$. Define
the map $f:\A\ra\A'$ by $f(t^i)=m^i$. Since the $\lambda^i$ are a
basis for $\m$, $f$ is completely determined by this requirement.
It is obvious that $f_*(\d+\mu_i\tns t^i)=\d+\mu_i\tns m^i$. It
is also evident that $f$ is unique. Thus $\d+\mu$ is a universal
infinitesimal deformation of $\d$.

Finally, suppose that \A\ is a \k-algebra which may not be
infinitesimal. The algebra $\A/\m^2$ is infinitesimal, with
maximal ideal $\m/\m^2$. Let $\tau:\A\ra\A/\m$ be the natural
projection, and suppose that $\td=d+\delta$ is an \A-deformation
of $d$. Then $d+\tau_*(\delta)$ is an infinitesimal deformation,
and thus determines an element $T(\delta)=[\tau_*(\delta)]$ in
$H(\L)\htns \m/\m^2$, which we call the differential of the
deformation  $\td$.

Suppose that \m\ is graded, and define the
\emph{tangent space} of \A\ by $\tanA=(\m/\m^2)^*$. Note that
$\tanA$ is complete. When
$\HL$ is complete and $\m$ is of finite type the differential can
viewed as a continuous map $T(\tilde\delta):\tanA\ra\HL$.

In the case of an infinitesimal deformation, we can express
$\delta=\delta_i\tns m^i$, where $\delta_i\in\ZL$, so that its
differential can be expressed as
$T(\td)=[\delta]=[\delta_i]\tns\bar m_i$. If the deformation is
not infinitesimal then we can still express
$T(\td)=[\delta_i\tns\bar m_i]$, in terms of a decomposition
$\delta=\delta_i\tns m_i$, but the expression $[\delta_i]\tns
\bar m_i$ may not make sense because $\delta\not\in\ZL\tns\m$.

Our main objective is to extend the universal infinitesimal
deformation to a miniversal deformation. For brevity, let us
denote $\H=(\Pi\HL)^*$.  Then the algebra $\A=\k\bigoplus\H$
corresponding to the universal infinitesimal deformation can be
expressed in the form $\A=\k[[\H]]/(\H)^2$. If we have any
$\A'$-deformation, then there is a unique map from \A\ to
$\A'/(\m')^2$ which takes the universal infinitesimal deformation
to the induced infinitesimal deformation. We can lift the map
(non-uniquely) to a map from $\A$ to $\A'$, which induces a unique
algebra morphism $\k[\H]$ to $\A'$.  When $\A'$ is a formal
algebra, this morphism extends uniquely to a morphism
$\k[[\H]]\ra A'$.  A problem that arises is that in general, the
universal infinitesimal deformation does not extend to a
$\k[[\H]]$-deformation; instead, we will need to consider a
quotient algebra of $\k[[\H]]$.

The non-universality of the deformation is related to the
necessity to choose a lift of $\A'/(\m')^2$ to $\A'$. It is
precisely this point which gives rise to the fact that there is
no universal object in the category of formal deformations.  Note
that an induced homomorphism from a quotient of $\k[[\H]]$ to
$\A'$ is unique if it exists, and when $\A'$ is infinitesimal,
this homomorphism is unique, because in this case there is no
lift, and therefore no choice to make.  Thus such a quotient
algebra is a good candidate for a miniversal deformation.  In
order to construct the appropriate algebra, we will need to
analyze how to extend an algebra by a trivial module.  This
subject was studied extensively by Harrison \cite{harr}, and we
shall explore this matter in the next section.

\section{Harrison Cohomology}

In order to understand how to construct a miniversal deformation, we
shall have to consider commutative extensions of a commutative algebra
by a module. We shall assume that all structures below are filtered
spaces and all maps are continuous.

 \emph{An extension  \B\ of an algebra \A\ by an
\A-module $N$} is a \k-algebra  \B\ together with an exact
sequence of \k-modules
\begin{equation}
0\ra N\stackrel{i}{\ra}\B\stackrel{p}{\ra} \A\ra 0,
\end{equation}
where $p$ is an \k-algebra homomorphism, and the  \B-module
structure on $i(N)$ is given by the \A-module structure of $N$ by
$i(n)b=i(n(p(b))$. In particular, if we identify $N$ with its
image $i(N)$, then $N$ is an ideal in  \B\ satisfying $N^2=0$.
Let $\lambda:\B\ra N\oplus \A$ be a \k-module isomorphism such
that $\lambda(n)=(n,0)$ for $n\in N$, and such that $p=\pi_2
\circ\lambda$, where $\pi_2$ denotes the projection onto the
second component. Then $\lambda$ determines a product on $N\oplus
\A$ such  that $\lambda(b)\lambda(b')=\lambda(bb').$ We have
\begin{equation*}
(n,0)(n',0)=\lambda(n)\lambda(n')=\lambda(nn')=(0,0).
\end{equation*}
If $\lambda(b)=(0,a)$, then $p(b)=a$, so we have
\begin{equation*}
(n,0)(0,a)=\lambda(n)\lambda(b)=\lambda(nb)=\lambda(na)=(na,0).
\end{equation*}
Finally
\begin{equation*}
(0,a)(0,a')=(\ph(a,a'),aa')
\end{equation*}
for some \k-linear map $\ph:a\tns a\ra N$, which has degree 0, and
is graded symmetric; that is, $\ph(a,a')=\s{aa'}\ph(a',a)$.
\emph{Two extensions  \B\ and  $\B'$ of \A\ by $N$ are said to be
equivalent} if there is an \k-algebra isomorphism $f: \B\ra\B'$
such that the diagram below commutes.
$$
\begin{CD}
0@>>>N@>>>\B@>>>\A@>>>0\\
&&@| @VfVV @|\\
0@>>>N@>>>\B'@>>>\A@>>>0\\
\end{CD}
$$

An equivalence from  \B\ to  \B\ is said to be an
\emph{automorphism of  \B over \A}. We use a graded version of
Harrison cohomology to characterize the properties of extensions.

First, consider the map $\ph$ above. It is easy to show that
associativity of the multiplication is equivalent to the cocycle
condition
\begin{equation*}
a\ph(b,c)-\ph(ab,c)+\ph(a,bc)-\ph(a,b)c=0.
\end{equation*}
This condition immediately yields $\ph(1,a)=\ph(1,1)a$. Next,
note that if $\lambda(e)=(-\ph(1,1),1)$, then $e$ is the
multiplicative identity in \B\ because if $\lambda(b)=(n,a)$, then
\begin{multline*}
\lambda(eb)=(-\ph(1,1),1)(n,a)=
\\(-\ph(1,1)a+n+\ph(1,a),a)=(n,a)=\lambda(b).
\end{multline*}
If we define $\lambda'$ by
$\lambda'(b)=\lambda(b)+(\ph(1,1)b,0)$, then
 $\lambda'(n)=(n,0)$ and $\pi_2 \circ\lambda'=p$.
Furthermore, if the product in terms of the decomposition of
 \B\ determined by $\lambda$ is given by the cocycle $\ph'$,
then since $\lambda'(e)=(0,1)$, we have
\begin{equation*}
(0,1)=(0,1)(0,1)=(\ph'(1,1),1),
\end{equation*}
so $\ph'(1,1)=0$.  Thus, if $a=k+m$, $a'=k'+m'$ are elements of
$\A$ given in terms of the decomposition $\A=\k\oplus\m$, then
$\ph'(a,a')=\ph'(m,m')$, so $\ph'$ is completely determined by
its restriction to $\m\tns\m$. We shall call a cocycle $\ph$
satisfying $\ph(1,1)$ a reduced cocycle, and we have shown that
every extension $\B$ can be defined by a decomposition $B=N\oplus \A$,
where the product is given by a reduced cocycle.

In order to relate this to Harrison cohomology, define
$Ch^2(\A,N)$ to be the submodule of $\hom(\m^2,N)$ consisting of
symmetric maps, and define $d_2:Ch^2(\A,N)\ra\hom(\m^3,N)$ by

\begin{multline}
d_2\ph(m,m',m'')=\\
\s{m\ph}m\ph(m',m'')-\ph(mm',m'')+\ph(m,m'm'')- \ph(m,m')m''.
\end{multline}
Note that even though deformations are determined only by even
cocycles, we do not restrict our definitions to such elements,
hence the sign appears in this definition of the coboundary
operator. Also, $Ch^1(\A,N)=\hom(\m,N)$ is the space of Harrison
1-cochains, with $d_1:Ch^1(\A,N)\ra Ch^2(\A,N)$ given by
\begin{equation}
d_1\lambda(m,m')=\s{m\lambda}m\lambda(m')-\lambda(mm')+\lambda(m)m'.
\end{equation}
It is easily checked that $d_1\lambda$ is graded symmetric, and that
$d_1^2=0$. The condition $d_1\lambda=0$ is just the derivation
property, so
\begin{equation*}
Ha^1(\A,N)=\ker(d_1)=\der(\A,N).
\end{equation*}
Define
\begin{equation*}
Ha^2(\A)=\ker(d_2)/\im(d_1).
\end{equation*}

We show that the even part of $Ha^1(\A,N)$ classifies the
automorphisms of an extension, while the even part of
$Ha^2(\A,N)$ classifies the equivalence classes of extensions.
Note that the even part of $Ha^k(\A,N)$ is not determined by just
the even part of $Ha^k(\A,\k)$, which explains why we need to
consider all of the cohomology, even though only even elements
actually give deformations and automorphisms.

Next, note that we have already shown that an extension of \A\ by
$N$ is given by a cocycle $\ph$ in $Ch^2(\A,N)$. We show that
extensions $\B$ and $\B'$, given by cocycles $\ph$ and $\ph'$ are
equivalent precisely when they differ by a coboundary. In terms
of the decompositions of $\B$ and $\B'$ determined by the
cocycles, an equivalence is given by a map $f:\B\ra\B'$ satisfying
$f(n,0)=(n,0)$ and $f(0,a)=(\lambda(a),a)$. If $f$ is a
homomorphism, then the two lines below are equal
\begin{align*}
f((0,a)(0,a'))&=f(\ph(a,a'),aa')=(\lambda(aa')+\ph(a,a'),aa'),\\
f(0,a)f(0,a')&= (\lambda(a),a))(\lambda(a'),a')=
(\lambda(a)a'+a\lambda(a')+\ph'(a,a'),aa'),
\end{align*}
so that $\ph=\ph'+d\lambda$. Conversely, if $\ph=\ph'+d\lambda$
for some even $\lambda\in Ch^1(\A,N)$, then
$f_\lambda(n,a)=(n+\lambda(a),a)$ defines an equivalence. Thus
two cocycles are equivalent precisely when they differ by a
coboundary, and we see that the even part of $Ha^2(\A,N)$
classifies the extensions of \A\ by $N$. Applying the analysis
above to an automorphism of \B\ over \A, we see that it is
determined by an element $\lambda\in Ch^1(\A,N)$ satisfying the
condition $d\lambda=0$. In other words, it is a Harrison
1-cocycle. Since $Ha^1(\A,N)$ is precisely the space of cocycles
(there are no 1-coboundaries), we see that the automorphisms of
\B\ over \A\ are classified by $Ha^1(\A,N)$.

A morphism between an extension \B\ of \A\ by $N$ and an
extension $\B'$ of \A\ by $N'$ is given by a commutative diagram
$$
\begin{CD}
0@>>>N@>>>\B@>>>\A@>p>>0\\
&&@VgVV @VfVV @|\\
0@>>>N'@>>>\B'@>>>\A@>p'>>0\\
\end{CD},
$$
where $g:N\ra N'$ is an \A-module homomorphism and $f:\B\ra \B'$
is an \A\ algebra homomorphism. Given $g$, the homomorphism $f$,
if it exists, is determined up to an automorphism of $\B'$. To see
this, suppose that $f$ and $f'$ both satisfy the above
requirements. Then $f(b)-f'(b)=n'$ for some $n'\in N'$, and if
$p(b)=p(b')$, then $b-b'=n$ for some $n\in N$, so that
\begin{equation*}
f(b)-f'(b)-f(b')+f'(b')=f(b-b')-f'(b-b')=g(n)-g(n)=0.
\end{equation*}
Thus we can define an even $\lambda:\A\ra N'$ such that
$f(b)-f'(b)=\lambda(p(b))$. Suppose that $p(b)=a$ and $p(b')=a'$.
Then we see that
\begin{multline*}
\lambda(aa')=f(bb')-f'(bb')=
f(b)f(b')-f'(b)f'(b')=\\
f(b)(f(b')-f'(b'))+(f(b)-f'(b))f'(b')=\\
f(b)\lambda(a')+\lambda(a)f(b')= a\lambda(a')+\lambda(a)a',
\end{multline*}
which shows that $\lambda$ is a derivation. Now recall that the
automorphism $f_\lambda$ is given by
$f_\lambda(b')=b'+\lambda(p'(b'))$. Finally, note that
\begin{equation*}
f_\lambda\circ f'(b)=f'(b)+\lambda(p'(f'(b)))=
f'(b)+\lambda(p(b))=f(b).
\end{equation*}
This shows that $f$ is uniquely defined up to an automorphism of
$\B$, so that the equivalence classes of such mappings are
determined by the map $g$.

Of course, we still have to determine when such a map $f$ exists.
Note that $g$ induces a map $g_*:Ch^n(\A,N)\ra Ch^n(\A,N')$, which
commutes with the differential, and so induces a map
$g_*:Ha^k(\A,N)\ra Ha^k(\A,N')$. Then we claim that $f$ exists
precisely when $g_*(\ph)$ is equivalent to $\ph'$, where $\B$ and
$\B'$ are given by the 2-cocycles $\ph$ and $\ph'$ respectively.
To see this, express $\B=N\oplus\A$ and $\B'=N'\oplus \A$, with
multiplication given by the cocycles. Suppose that $f:\B\ra \B'$
exists. Then $f(n,a)=(g(n)+\lambda(a),a)$, and the homomorphism
property shows that
\begin{multline*}
(g\circ\ph(a,a')+\lambda(aa'),aa')= f((0,a)(0,a'))=
f(0,a)f(0,a')=\\
(\lambda(a),a)(\lambda(a'),a')=
(a\lambda(a')+\lambda(a)a'+\ph'(a,a'),aa'),
\end{multline*}
so that $g_*\ph=\ph'+d\lambda$, and they are equivalent.
Conversely, if this relation holds, then
$f(n,a)=(g(n)+\lambda(a),a)$ is a homomorphism satisfying the
requirements. Thus it makes sense to define a morphism of
extensions as a map $g:N\ra N'$ which extends to a map $f:\B\ra
\B'$, having the required properties, since this only depends on
the equivalence classes of the extensions $\B$ and $\B'$. Note
that in general, there may be many such maps $g$.

When \A\ is an augmented algebra with augmentation ideal \m, then
an \emph{infinitesimal extension of \A\ by a \k-vector space $N$}
is one in which the \A-module structure of $N$ is given by the
augmentation, so that in particular, $N\cdot\m=0$.

Now, let us suppose that \m\ is a \gr\ space. Then it is
natural, when $N$ is also \gr, to let $Ch^*(\A,N)$ be the
cochains respecting the grading.  When \m\ is of finite type, there is an
initial object in the category of extensions of \A\ by \gr\
\A-modules $N$, which will play a role in the construction of the
miniversal deformation of $L$. The inclusion
$\hom(\m^k,N)\subseteq {N\htns (\m^k)}^*$ gives rise to a map
\begin{equation*}
Ch(\A,N)\ra N\htns Ch^k(\A,\k)
\end{equation*}
which commutes with the differential when $N$ is infinitesimal, and so induces a map
\begin{equation*}
Ha^k(\A,N)\ra N\htns Ha^k(\A,\k).
\end{equation*}

When $N$ is of finite type, these maps are isomorphisms.  Note that
$Ha^k(\A,\k)$ is a complete filtered space, which can be seen as
follows. 
First, note that $Ch^k(\A,\k)$ is a complete subspace of
${(\m^k)}^*$. Define a map $b:\m^{k+1}\ra \m^k$ by
\begin{equation*}
b(m_1\mcom m_{k+1})=\sum_{i=1}^{k}(-1)^i(m_1\mcom m_im_{i+1},m_{i+2}\mcom
m_{k+1}).
\end{equation*}
It is easy to see that the Harrison coboundary operator
the dual of the map above, so is 
continuous.  Thus the space of coboundaries is closed, and $Ha^k(\A,\k)$
is naturally a complete filtered space.

Let $M=Ha^2(\A,\k)^*$, equipped with the \A-module structure
determined by the augmentation,  and choose some order preserving
$\mu:Ha^2(\A,\k)\ra Ch^2(\A,\k)$ such that
$\mu(\bar\ph)\in\bar\ph$. Then $\mu^*:\m\tns\m\ra M$, given by
\begin{equation}
\mu^*(m,m')(\bar\ph)=\s{(m+m')\ph}\mu(\bar\ph)(m,m')
\end{equation}
is an even, order preserving 2 cochain. It is a 2-cocycle
because
\begin{multline*}
d\mu^*(m,m',m'')(\bar\ph)=\\
m\mu^*(m,m')(\bar\ph)-
\mu^*(mm',m'')(\bar\ph)+\mu^*(m,m'm'')(\bar\ph)+\mu^*(mm')m''(\bar\ph)\\
=\s{(m'+m'')\ph}m\mu(\bar\ph)(m',m'')-
\s{(m+m'+m'')\ph}\mu(\bar\ph)(mm',m'')+\\
\s{(m+m'+m'')\ph}\mu(\bar\ph)(m,m'm'')-
\s{(m+m')\ph+m''\ph}\mu(\bar\ph)(m,m')m''=\\
\s{(m+m'+m'')\ph}d\mu(\bar\ph)(m,m',m'')=0.
\end{multline*}
Thus $\mu^*$ determines an extension $0\ra M\ra \mbox{$\mathcal C$}
\ra \A\ra 0$ of
\A\ by $M$. This extension does not depend, up to equivalence, on
the choice of $\mu$. To see this suppose that $\mu'$ is another
choice, so that $\mu(\bar\ph)-\mu'(\bar\ph)=d\psi$ for some
$\psi\in Ch^1(\A)$. Note that $|\ph|=|\psi|$. We can define
$\lambda:\m\ra M$ by
\begin{equation*}
\lambda(m)(\bar\ph)=\s{m\ph}\psi(m).
\end{equation*}

Then
\begin{multline*}
d\lambda(m,m')(\bar\ph)=
m\lambda(m')(\bar\ph)-\lambda(mm')(\bar\ph)+\lambda(m)m'(\bar\ph)=\\
\s{m'\ph}m\psi(m')-\s{(m+m')\ph}\psi(mm')+\s{(m+m')\ph}\psi(m)m'=\\
\s{(m+m')\ph}d\psi(m,m')= (\mu^*-{\mu'}^*)(m,m')(\bar\ph).
\end{multline*}

Now we show that this extension, when \m\ is  of finite type, is
universal in the set of infinitesimal extensions. Let $0\ra N\ra
\B\ra \A\ra 0$ be an extension given by some (even) $\bar\ph\in
Ha^2(\A,N)$. Using the inclusion $Ha^2(\A,N)\subseteq N\tns
{Ha^2(\A,k)}^*$, we can express $\bar\ph=n_i\tns\bar\ph^i$.
Define the map $g:M\ra N$ by
$g(\eta)=\s{\eta\ph^i}n_i\eta(\bar\ph^i)$. Let
$\ph^i=\mu(\bar\ph^i)$ and define $\ph=n_i\tns \ph^i$. Then $\ph$
is a cocycle representing the cohomology class $\bar\ph$, and we
may assume that the decomposition of $\B=\A\oplus N$ is given by
the cocycle $\ph$.

Let $\ph\in\bar\ph$ be chosen so that $\ph^i=\mu^*(\bar\ph^i)$.
Then in terms of the decompositions
of $\B$ given by $\ph$ and ${\mbox{$\mathcal C$}}$ given by $\mu^*$, we have

\begin{multline*}
g_*(\mu^*)(m,m')= (g\circ\mu^*)(m,m')=
g(\mu^*(m,m'))=\\
\s{n_i\mu^*(m,m')}n_i\mu^*(m,m')(\bar\ph^i)=\\
\s{n_i\mu^*(m,m')+(m+m')\ph^i}n_i\ph^i(m,m')= \ph(m,m'),
\end{multline*}
because the signs cancel owing to the fact that $\ph$ and $\mu^*$
are even maps. Any other $g:M\ra N$ would determine a non
equivalent cocycle, so the morphism is unique. The extension of
\A\ by $M$ will be called the universal infinitesimal extension
of \A. An extension of \A\ by a module $N$ is called
\emph{essential} when the cocycle $\ph:\m\tns\m\ra N$ is
surjective. It is useful to note that the map $\mu^*:\m\tns\m\ra
M$ is surjective, so that it is an essential extension of \A.

Finally, we introduce the notion of a \gr\ algebra \A, and
the space of order preserving cochains.  In the definition of a
\gr\ algebra, we require that $\A^k\cdot \A^l\subseteq
\A^{k+l}$, and for augmented \gr\ algebras, we require that
$\m=\A^1$, so that nonzero elements in the field \k\ have order
zero. A typical example is a polynomial algebra, where the
generators are taken to be elements of nonzero degrees. A
\gr\ module $N$ over \A\ is required to satisfy $N^k\cdot
\A^l\subseteq N^{k+l}$. We are interested in classifying the
\gr\ extensions of \A\ by a \gr\ module $N$.

Since the augmentation ideal in the extended algebra is simply
$\m'=\m\oplus N$, it is necessary that $N=N^1$. Let $\ph$ be the
cocycle determining the extension. If $(0,m)$ and $(0,m')$ are
elements in the extended algebra, with $m$ $m'\in \m$, then since
$\ph(m,m')=(0,m)\cdot(0,m')$, it follows that $\ph(m,m')\in
N^{k+l}$ if $m\in M^k$ and $m'\in M^l$. Let us say that a map
$f:F\ra G$ between two \gr\ spaces is order preserving if
$f(F^n)\subseteq G^n$ for all $n$. If we consider a 2-cochain
$\ph$ to be a map $\ph:\m\tns\m\ra N$, our requirement is simply
that $\ph$ is order preserving. If $\ph\in
Ch^k(\A,N)$ is order preserving, it is easy to see that $d\ph$ is
also order preserving.

Now let us consider the universal infinitesimal extension of a
\gr\ algebra. Then since $\mu^*$ is both order preserving and
surjective, it follows that the universal extension is also a
\gr\ algebra, and is thus a universal object in the category
of extensions of \gr\ algebras.

Finally, let us suppose that $0\ra M\ra \B\ra \A\ra 0$ is the
universal infinitesimal extension of \A, and that $f:\A\ra\A'$ is
an algebra homomorphism.  Let $\B'$ be an infinitesimal extension
of $\A'$ by an $\A'$-module $N$. Then there is a unique extension
of the homomorphism $f$ to a homomorphism $f':\B\ra\B'$. To see
this, note that $N$ inherits an \A-module structure through $f$,
and moreover, if $\B'$ is decomposed in the form $\B'=\A'\oplus N$
using the cocycle $\ph'$, then $\ph=\ph'\circ(f\tns f)$
determines a cocycle in $Ha^2(\A,N)$. Thus we obtain an
infinitesimal extension $\B''$ of \A\ by $N$, and an obvious
homomorphism $\A\oplus N\ra \A'\oplus N$ extending $f$.  But there
is a homomorphism from the universal extension of \A\ to this
extension of \A, and composition of the two homomorphisms yields
the desired map. The commutative diagram below summarizes this
construction:

$$
\begin{CD}
0@>>>M@>>>\B@>>>\A@>>>0\\
&&@VVV @VVV @|\\
0@>>>N@>>>\B''@>>>\A@>>>0\\
&&@|@Vf'VV @VfVV\\
0@>>>N@>>>\B'@>>>\A'@>>>0\\
\end{CD}
$$
When $f$ is surjective and $N$ is an essential extension of
$\A'$, then it can be seen that $f'$ is also surjective. Note
that when the algebras are \gr, and $f$ is order preserving,
then we obtain an order preserving extension of this homomorphism.

An example of an infinitesimal extension which will be important
to us later arises when \A\ is a formal algebra, with maximal
ideal \m, and we let $\A_k=\A/\m^k$.  Then $N_k=\m^k/\m^{k+1}$ is
naturally an infinitesimal \A-module, and we have an exact
sequence
\begin{equation*}
0\ra N_k\ra \A_{k+1}\ra \A_k\ra 0,
\end{equation*}
expressing $\A_{k+1}$ as an essential, infinitesimal extension of
$\A_k$ by $N_k$, when $k\ge 1$. Now consider the formal algebra
$\A'=\k[[X]]$, where $X=\m/\m^2$, and its corresponding quotient
algebras $\A'_k=\A'/[X]^k$ ($[X]$ is the maximal ideal in $\A'$).
By the universal properties of the algebra $\A'$, one sees that
the extension of $\A'_k$ by $N'_k=[X]^k/[X]^{k+1}$ is the
universal infinitesimal extension of $\A'_k$. Moreover,
$\A_2=\k[[X]]/[X]^2=\A'_2$, so that the identity extends to a
homomorphism $f_3:\A'_3\ra A_3$, which is surjective, because
both extensions are essential.  Continuing on, we obtain a
sequence of surjective homomorphisms $f_k:\A'_k\ra \A_k$,
compatible with the projections between these algebras. It
follows that $\A$ is a quotient algebra of $\k[[X]]$, and
moreover, we see that if $\A_k=\k[[X]]/I_k$ for some ideal $I_k$,
then $[X]^k\subseteq I_k\subseteq [X]^2$.

A stronger characterization of the ideals $I_k$ above can be
obtained from some results due to Harrison, (propositions 5.1-2 in
\cite{ff2}).
\begin{thm}
Let $\A=\k[x_1\mcom x_n]$ be a polynomial algebra, and $\m$ be
the ideal generated by $x_1\mcom x_n$. If $I$ is an ideal in \A\
contained in $\m^2$, then
\begin{equation*}
Ha^2(\A/I,\k)=(I/(\m\cdot I))^*.
\end{equation*}

Furthermore, the universal infinitesimal extension of $\A/I$
is given by the exact sequence
\begin{equation*}
0\ra I/(\m\cdot I)\ra \A/(\m\cdot I)\ra \A/I\ra 0.
\end{equation*}
\end{thm}
The generalization of this result to the case where $\A=\k[[X]]$
for some \gr, finite type \k-vector space $X$ is straightforward.

Finally, let us suppose that $0\ra N\ra\B\ra\A\ra0$ is an
infinitesimal extension of \A\ and that $N'$ is a closed subspace
of $N$.  Then $N'$ is an ideal in \B, and we obtain an extension
$0\ra N/N'\ra\B/N'\ra\A\ra0$. When the extension given by $N$ is
essential, so is the extension by $N/N'$.

\section{Obstructions to extensions}
Let $\td=\d+\delta$ be a deformation of $\d$ with base \A\, and
suppose that $0\ra N\ra\B\ra\A\ra0$ is an infinitesimal extension
of \A. If we extend the base of the deformation to \B, the
Maurer-Cartan formula (\refeq{mc formula}) will not hold in
general, but instead we obtain that
\begin{equation}
\gamma=D(\delta)+\tfrac12\br\delta\delta\in \L\htns N.
\end{equation}
Moreover, $\gamma$ is a cocycle in $\L\htns N$, which can be seen
as follows. First note that for an odd element in a \zt-graded
Lie algebra, while graded antisymmetry does not force the bracket
of the element with itself to vanish, nevertheless, triple
brackets do vanish, \ie $\br\delta{\br\delta\delta}=0$. Next,
note that $\br{D(\delta)+\br\delta\delta}\delta=0$ because the
first term in the bracket is in $\L\htns N$ and the second lies in
$\L\htns\m$, and $N\m=0$. Using these facts, we obtain
\begin{equation*}
D(D(\delta)+\tfrac12\br\delta\delta)= \tfrac12D\br\delta\delta=
\br{D(\delta)}\delta=
\br{D(\delta)+\tfrac12\br\delta\delta}\delta=0.
\end{equation*}
In order to extend $\td$ to $\LB$, we still can add a term
$\beta\in\L\htns N$, and it is easy to see that  $\td+\beta$ is a
codifferential in \LB\ precisely when
\begin{equation*}
D(\beta)=-\gamma.
\end{equation*}
Thus the cohomology class $\bar\gamma$ in $H(\L)\htns N$
determines an obstruction to extending $\td$ to a deformation of
$\d$ with base \B. Note that in the case of a Lie algebra, the
obstruction actually lies in $H^3(\L)\htns N$, but in general, we
can only say that it is an even element in $H(\L)\htns N$.

To see that the obstruction depends only on the element $f\in
Ha^2(\A)$ which determines the extension, suppose that the
extension $\B$ is given explicitly by the Harrison cocycle
$\ph:\m^2\ra N$. Writing $\delta=\delta_i\tns m^i$, it is easily
seen that
$\gamma=\s{\delta_jm^i}\br{\delta_i}{\delta_j}\tns\ph(m^i,m^j)$.
If $\gamma'$ is the cocycle determined by the extension $\B'$
given by $\ph'$, and we express $\ph-\ph'=\delta(\lambda)$, for
some $\lambda:\m\ra N$, then it follows that
\begin{equation*}
\gamma-\gamma'=\s{\delta_jm^i}\br{\delta_i}{\delta_j}\tns\lambda(m^im^j)=
D(\delta_i\tns\lambda(m^i)),
\end{equation*}
and is thus a coboundary. Thus we have shown that the obstructions
to extending the deformation to an extension of \A\ by $N$
determines a map $\O:Ha^2(\A,N)\ra H(\L)\htns N$.

Supposing that the obstruction vanishes, the element $\beta$
constructed above is determined only up to a cocycle $\psi\in
\ZL\htns N$. Moreover, adding an odd coboundary in $\BL\htns N$
produces an equivalent deformation, as in the case of
infinitesimal deformations, owing to the fact that $N^2=0$. To
see this explicitly, let $\lambda\in\L\htns N$ be an even
coderivation. Then $f=\id+\lambda$ is a coderivation of $\SW\tns
\B$ (fixing $\SWA$) and $f\inv=\id-\lambda$, due to the fact
that $N$ is infinitesimal. Then
$f_*(d+\delta+\beta)=d+\delta+\beta+D(\lambda)$. Moreover, if
$d+\delta+\beta$ and $d+\delta+\beta'$ are equivalent extensions,
then adding a cocycle $\psi\in\ZL\htns N$ to each of them produces
an equivalent deformation. Thus we obtain a transitive action of
$\HL\htns N$ on the equivalence classes of extensions of the
deformation $\td$ to $\L\htns \B$. Note that we have not claimed
that nonequivalent cocycles give rise to nonequivalent
\B-deformations of $d$, because we only considered
equivalences arising from coderivations of $\SW\tns \B$ fixing
$\SWA$.

Now the set of automorphisms of \B\ over \A\ also acts on the
equivalence classes of extensions of the deformation $\td$ to
$\L\htns\B$.  Let $f$ be an automorphism of \B\ over \A. Then in
terms of a decomposition $\B=\A\oplus N$, with
$f(n,a)=(n+\lambda(a),a)$ for some Harrison 1-cocycle $\lambda$,
we have $f_*(d+\delta+\beta)=d+\delta+\beta +\lambda_*(\delta)$.
Note that
\begin{equation*}
D(\lambda_*(\delta))=\lambda_*(D(\delta))
=\lambda_*(-D(\beta)-\frac12[\delta,\delta])
=\lambda_*(-\frac12[\delta,\delta])=0,
\end{equation*}
since $D(\beta)\in\L\htns N$, $[\delta,\delta]\in\L\htns\m^2$ and
$\lambda$, being a cocycle, vanishes on $\m^2$. Thus
$\beta+\lambda_*(\delta)$ determines another extension. If
$\beta$ and $\beta'$ are two equivalent extensions, then
$\beta+\lambda_*(\delta)$ and $\beta' +\lambda_*(\delta)$ are also
equivalent, so we see that $Ha^1(\A,N)$ acts on the equivalence
classes of extensions. We would like to show that this action is
transitive, because that would show that up to an automorphism of
\B, the extension of $\td$ to \B\ is unique up to equivalence. To
do this, we relate the transitive action of $H(\L)\htns N$ to the
action of $Ha^1(\A,N)$. Since $\lambda_*(\delta)\in Z(\L)\htns
N$, it determines an element in $H(\L)\htns N$. We investigate
when this map is surjective.

The action of  $\lambda\in Ha^1(\A,N)$ on $\delta$ can be thought
of as an action of the differential $T(\delta)\in
H(\L)\htns\m/\m^2$ on $Ha^1(\A,N)$, since $\lambda$ vanishes on
$m^2$, and so acts on $\m/\m^2$. Let us say that the differential
is surjective if we can express $T(\delta)=\delta_i\htns m^i$,
where $\delta_i$ is an increasing spanning subset of $H(\L)$, and
$m^i$ is an increasing basis of $\m/\m^2$. When $\m/\m^2$ is a
dual space of a filtered space $F$, then this is equivalent to
saying that the associated map in $\hom(F,H(\L))$ is surjective.
Then we claim that if $N$ is an infinitesimal extension of
$d+\delta$, and $T(\delta)$ is surjective, then the differential
induces a surjective map $Ha^1(\A,N)\ra H(\L)\tns N$. To see
this, note that since $\{\delta_i\}$ spans $H(\L)$, any element
$H(\L)\htns N$ can be written in the form $\delta_i\htns n_i$.
Then define $\lambda\in Ha^1(\A,N)$ by $\lambda(m^i)=n^i$ (using
the natural identification $Ha^1(\A,N)=\hom(\m/\m^2,N)$).  Then
it is clear that $\lambda_*(\delta)= \delta_i\htns n^i$.

Our only application of the above result will be to a special
sequence of extensions of the universal infinitesimal extension
which satisfies $\m/\m^2=(\Pi(H(\L))^*=\H$, as long as at each step
of the way, the extensions are infinitesimal and essential.  In
this case, the differential is surjective, since it is given by
the identity map $\Pi(H(\L))\ra H(\L)$.

One final result about the obstruction will prove useful in our
construction of the miniversal deformation.
\begin{thm}
Suppose that $\td=d+\delta$ is an \A-deformation of the codifferential
$d$ on $L$, and that
$$
\begin{CD}
0@>>>M@>>>\B@>>>\A@>>>0\\
&&@VgVV @VfVV @|\\
0@>>>N@>>>\B'@>>>\A@>>>0\\
\end{CD}
$$
is a diagram representing a morphism of extensions.  If
$\bar\gamma\in H(\L)\htns M$ is the obstruction to the extension
of $\td$ to $\L\htns\B$, then $\overline{g_*(\gamma)}$ is the
obstruction to the extension of $\td$ to $\L\htns\B'$.
\end{thm}
\begin{proof}
If $\delta=\delta_i\htns m^i$, then
$\gamma=\s{\delta_jm^i}[\delta_i,\delta_j]\htns \ph(m^im^j)$,
where $\ph$ is the cocycle determining the extension to $\B$.
Similarly, if $\bar\gamma'$ is the obstruction to extending the
deformation to $\B'$, we can express
$\gamma'=\s{\delta_jm^i}[\delta_i,\delta_j]\htns \ph'(m^im^j)$,
where $\ph'$ determines the extension $\B'$.  Since
$g(\ph(m^im^j))=\ph'(m^im^j)$, the result follows immediately.

\end{proof}

\section{Construction of the Miniversal Deformation}

We have assembled all of the tools for constructing the
miniversal deformation in the previous two sections, so all that
remains here is to tie the ideas together. Suppose that $L$ is
a differential graded Lie algebra, $d$ is a
codifferential on $\L$, $\A'$ is a filtered formal algebra with
maximal ideal $\m'$, and that $\td=d+\delta'$ is a formal
deformation of $d$.  Let $\A'_k=\A/(\m')^k$, and
$\td'_k=d+\delta'_k$ be the induced deformation of $\td$ on
$\A'_k$. We will construct a formal algebra $\A$ with maximal
ideal $\m$, formal deformation $\td=d+\delta$ with induced
deformation $\td_k$ on $\L\htns \A_k$, and homomorphisms
$f_k:\A_k\ra \A'_k$. Here $\A_k=\A/\m^k$ are such that the
diagram of extensions commutes
$$
\begin{CD}
0@>>>N_k@>>>\A_{k+1}@>>>\A_k@>>>0\\
&&@Vg_kVV @Vf_{k+1}VV @Vf_kVV\\
0@>>>N'_k@>>>\A'_{k+1}@>>>\A'_k@>>>0\\
\end{CD},
$$
where $N'_k=(\m')^k/(\m')^{k+1}$, and in addition,
$(f_k)_*(\td_k)=\td'_k$, so that the induced formal homomorphism
$f:\A\ra\A'$ satisfies $f_*(\td)=\equiv\td'$.

For $k=1$, $\A_1=\A'_1=\k$, and $\td_1=\td'_1=d$. Next, for $k=2$, we
have $N_1=\H$, and the extension $\td_1=d+\delta_1$ is the
universal infinitesimal extension of $d$, which comes equipped
with a homomorphism $g_2:N_1\ra N'_1$, since the extension of
$\A'_1$ by $N'_1$ is infinitesimal. Now suppose that we have
constructed $\A_k$ and the map $f_k$ satisfying the requirements.
Then consider the universal extension $M$ of $\A_k$ given by
Harrison cohomology. Consider the induced extension of
$\A_k$ determined by extending it by $N'_k$ using the
homomorphism $f_k$ to define the extension. Then we claim that
the obstruction to extending $\td_k$ to the extension of $\A_k$
by $N'_k$ is the same as the obstruction to extending the
deformation $\td'_k$ to $\A'_{k+1}$, and therefore vanishes. But
then consider the morphism of extensions

$$
\begin{CD}
0@>>>M@>>>\A_k\oplus M@>>>\A_k@>>>0\\
&&@VgVV @VfVV @Vf_kVV\\
0@>>>N'_k@>>>\A_k'\oplus N'_k@>>>\A_k'@>>>0\\
\end{CD},
$$
and suppose that $\bar\gamma=\delta_i\htns n^i$ is the
obstruction to extending $\td_k$ to $\A_k\oplus M$. Then
$g_*(\bar\gamma)=0$, so that in particular $g(n^i)=0$ for all
$i$.  Let $M'$ be the closed subspace of $M$ generated by the
elements $n_i$. The map $g$ factors through the quotient
$N_k=M/M'$,and the obstruction to extending the deformation to
the extension $0\ra N_k\ra \A_k\oplus M/M'\ra \A_k\ra 0$ clearly
vanishes. Let $\A_{k+1}=\A_k\oplus N_k$. Since the extension of
$\A_k$ by $N'_k$ is essential, we know that any two extensions of
$\td_k$ to an $\A_k\oplus N'_k$ are equivalent up to an
automorphism. Take any extension $\td_{k+1}$ of $\td_k$ to
$\L\htns\A_{k+1}$ (again such an extension is determined up to
equivalence and automorphism). After applying an automorphism
$\eta$ of $\A_k\oplus N'$, we obtain the deformation
$\eta_*(f_*(\td_{k+1}))$, which is equivalent to a deformation
projecting to $\td'_{k+1}$ by the natural map $\A_k\oplus
N'\ra \A'_{k+1}$. Let $f_{k+1}=\eta\circ f$, then
$(f_{k+1})_*(\td_{k+1})\equiv \td'_{k+1}$.

We summarize our procedure in the following theorem.

\begin{thm} Suppose that $\L$ is a complete filtered \zt-graded
Lie algebra of finite type with codifferential $d$.  Then there
is a miniversal deformation $\td=d+\delta$ of $d$.
\end{thm}

The fact that the deformation constructed above is miniversal,
rather than just versal, is immediate from the manner in which it
is constructed by extensions of the universal infinitesimal
deformation.  Also, from some remarks made before, one can see that
if \A\ is the base of the miniversal deformation, then it is a
quotient of $\k[[\H]]$ by an ideal contained in $[\H]^2$, where
$\H=(\Pi(H(\L))^*$.  Furthermore, this ideal has an increasing
sequence of generators.

The proof of the existence of the miniversal deformation can be
considered constructive. In fact, at each stage one obtains an
obstruction $\gamma_k$, and this gives rise to some elements in
$[\H]^k$ that need to be set equal to zero.  Then one must solve the
problem $D(\beta)=\gamma$, which is not such an easy problem to solve
in practice.
Generalized Massey products (see \cite{fl}, \cite{pw}) play a role in
  this solution.

In a later paper, the authors will present some examples of the
construction of the miniversal deformation of \linf\ algebras.
It should be mentioned that even in the case of Lie algebras,
the construction of the miniversal deformation may  not be easy.  In
\cite{ff2}, \cite{fp}
 miniversal deformations of the infinite dimensional vector field
Lie algebras $L_1$ and $L_2$ are constructed, and even though the
cohomology is finite dimensional, the constructions are not simple.
Thus, the authors feel that carrying out detailed computations here,
while somewhat enlightening in the simplest cases, would overwhelm the
simplicity of the general results we have presented in this paper.
\bigbreak

\noindent\emph{Acknowledgements}

\noindent The authors wish to thank Dmitry Fuchs for helpful
discussions.
The first author thanks the Erwin Schr\"odinger Institute for
Mathematical Physics Vienna for partial support. The second author
thanks the E\"otv\"os Lor\'and University Budapest for hospitality.

\bibliographystyle{amsplain}
\providecommand{\bysame}{\leavevmode\hbox
to3em{\hrulefill}\thinspace}

\end{document}